\documentclass[11pt]{amsart}

\usepackage{amsrefs}
\usepackage[all]{xy}
\usepackage{syntonly}
\usepackage{hyperref}
\usepackage{amsfonts}
\usepackage{amssymb}
\usepackage{amsmath}
\usepackage{amsthm}
\usepackage[inline]{enumitem}
\usepackage{tikz}
\usepackage{graphics}
\usepackage{graphicx}
\usepackage{color}
\usepackage{comment}
\usepackage{caption,lipsum}

\usepackage{lineno}

\usepackage{multirow}

\newtheorem{theorem}{Theorem}[section]
\newtheorem{corollary}[theorem]{Corollary}

\newtheorem{definition}[theorem]{Definition}
\newtheorem{lemma}[theorem]{Lemma}

\newtheorem{remark}[theorem]{Remark}

\specialcomment{com}{ \color{red}\Large }{ \normalsize \color{black}} 
\specialcomment{oldcom}{ \color{brown} }{\color{black}} 

\excludecomment{oldcom} 

\usepackage{stmaryrd}

\begin{document}
\title[]{The $\Pi^1_1  \! \! \downarrow$ L\"owenheim-Skolem-Tarski property of Stationary Logic}

\author{Sean Cox}
\email{scox9@vcu.edu}
\address{
Department of Mathematics and Applied Mathematics \\
Virginia Commonwealth University \\
1015 Floyd Avenue \\
Richmond, Virginia 23284, USA 
}

\thanks{The author thanks Hiroshi Sakai (JSPS Kakenhi Grant Number 18K03397) for travel support to the RIMS Set Theory Workshop 2019, and also gratefully acknowledges support from Simons Foundation grant 318467.}


\begin{abstract}
Fuchino-Maschio-Sakai~\cite{FuchinoEtAl_DRP_LST} proved that the L\"owenheim-Skolem-Tarski (LST) property of Stationary Logic is equivalent to the Diagonal Reflection Principle on internally club sets ($\text{DRP}_{\text{IC}}$) introduced in \cite{DRP}.  We prove that the restriction of the LST property to (downward) reflection of $\Pi^1_1$ formulas, which we call the $\Pi^1_1 \! \! \downarrow$-LST property, is equivalent to the \emph{internal} version of DRP from \cite{Cox_RP_IS}.   Combined with results from \cite{Cox_RP_IS}, this shows that the $\Pi^1_1 \! \! \downarrow$-LST Property for Stationary Logic is strictly weaker than the full LST Property for Stationary Logic, though if CH holds they are equivalent.  
\end{abstract}

\maketitle

\tableofcontents

\section{Introduction}

\textbf{Stationary Logic} is a relatively well-behaved fragment of Second Order Logic introduced by Shelah~\cite{MR376334}, and first investigated in detail by Barwise et al~\cite{MR486629}.  Stationary Logic augments first order logic by introducing a new second order quantifier \emph{stat}; we typically interpret ``$\text{stat} Z \ \phi(Z,\dots)$" to mean that there are stationarily many countable $Z$ such that $\phi(Z,\dots)$ holds.\footnote{Other interpretations, e.g.\ for uncountable $Z$, or for filters other than the club filter, are often considered too.}  The quantifier \emph{aa} stands for ``almost all" or ``for club many"; so 
\[
\text{aa} Z \ \phi(Z,\dots)
\]
is an abbreviation for
\[
\neg \ \text{stat} Z \ \neg \ \phi(Z,\dots). 
\]
Section \ref{sec_prelims} provides more details.

By \emph{structure} we will always mean a first order structure in a countable signature.  The question of whether every structure has a ``small" elementary substructure in Stationary Logic was raised already in \cite{MR486629}.  One cannot hope to always get countable elementary substructures; e.g.\ if $\kappa$ is regular and uncountable, then $(\kappa,\in)$ satisfies ``$\in$ is a linear order and 
\[
\text{aa} Z  \ \exists x \ x \text{ is an upper bound of } Z",
\]
but no countable linear order can satisfy that sentence.  In a footnote in \cite{MR486629}, it was observed that even the statement 
\begin{equation}\label{eq_LST}
\begin{split}
& \text{``Every structure has an elementary (w.r.t.\ Stationary Logic)} \\
& \text{substructure of size } \le \omega_1 \text{"}
\end{split}  \tag{LST}
\end{equation}
 carries large cardinal consistency strength.\footnote{See Definition \ref{def_LST} for precisely what is meant by ``elementary substructure" in this context.}  The quoted statement above is now typically called the \emph{L\"owenheim-Skolem-Tarski} (LST) property of Stationary Logic.\footnote{The weaker assertion that every consistent theory (in Stationary Logic) has a model of size $\omega_1$, on the other hand, is a theorem of ZFC, as proven in \cite{MR486629}.}

Fuchino et al.\ recently proved that LST is equivalent to a version of the Diagonal Reflection Principle introduced in Cox~\cite{DRP}:
\begin{theorem}[Fuchino-Maschio-Sakai~\cite{FuchinoEtAl_DRP_LST}]\label{thm_LST_Fuchino}
LST is equivalent to the \textbf{Diagonal Reflection Principle on internally club sets} ($\text{DRP}_{\text{IC}}$).  
\end{theorem}

The purpose of the present note is to prove the following variant of Theorem \ref{thm_LST_Fuchino} involving $\Pi^1_1$ formulas in Stationary Logic (defined in Section \ref{sec_prelims} below) and the principle $\text{DRP}_{\text{internal}}$ from \cite{Cox_RP_IS}:
\begin{theorem}\label{thm_Cox_DRP}
The $\Pi^1_1 \! \! \! \downarrow$-LST property of Stationary Logic (see Definition \ref{def_LST}) is equivalent to the principle $\text{DRP}_{\text{internal}}$.
\end{theorem}

Cox~\cite{Cox_RP_IS} proved that $\text{DRP}_{\text{IC}}$ is strictly stronger than $\text{DRP}_{\text{internal}}$.  This was obtained by forcing over a model of a strong forcing axiom in a way that preserved $\text{DRP}_{\text{internal}}$ while killing $\text{DRP}_{\text{IC}}$ (in fact killing $\text{RP}_{\text{IC}}$; the argument owed much to Krueger~\cite{MR2674000}).  Furthermore, if CH holds, then $\text{DRP}_{\text{IC}}$ is equivalent to $\text{DRP}_{\text{internal}}$.  Combining those results with Theorem \ref{thm_Cox_DRP} immediately yields:
\begin{corollary}
The LST property of Stationary Logic is \textbf{strictly} stronger than the $\Pi^1_1 \! \! \downarrow$-LST property of Stationary Logic.

However, if the Continuum Hypothesis holds, they are equivalent.\footnote{One doesn't actually need the full continuum hypothesis for this equivalence to hold, but rather a variant of Shelah's Approachability Property, namely that the class of internally stationary sets is the same (mod NS) as the class of internally club sets.  See Cox~\cite{Cox_RP_IS} for more details.}
\end{corollary}

\begin{figure}
\caption{An arrow indicates an implication, an arrow with an X indicates a non-implication}  
\label{fig_LST}
\begin{displaymath}
\xymatrix{
\text{MM}^{++} \ar[dd] \ar[r] & \text{DRP}_{\text{IC}} \ar@{<->}[r] \ar[dd] &   \txt{Stationary Logic has\\ the LST property}  \\
& & \txt{(\emph{these 4 statements}\\\emph{are equivalent if CH holds})}  \\
\txt{Martin's\\Maximum}   \ar@/_0pc/!<0ex,0ex>;[r]!<0ex,0ex>|{\textcircled{\scriptsize \textbf{X}}} & \text{DRP}_{\text{internal}} \ar@{<->}[r]     \ar@/_0pc/!<2ex,0ex>;[uu]!<2ex,0ex>|{\textcircled{\scriptsize \textbf{X}}}  &   \txt{Stationary Logic has\\ the $\Pi^1_1 \! \! \downarrow$ LST property} \\
}
\end{displaymath}
\end{figure}
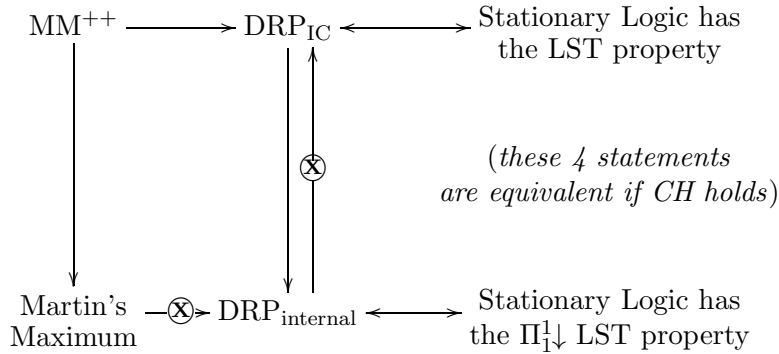

We note that while the technical strengthening $\text{MM}^{++}$ of Martin's Maximum implies $\text{DRP}_{\text{IC}}$ (see \cite{DRP}), recent work of Cox-Sakai~\cite{Cox_Sakai_DRP} shows that Martin's Maximum alone does not imply even the weakest version of DRP.  Figure \ref{fig_LST} summarizes the relevant implications and non-implications discussed in this introduction.

Section \ref{sec_prelims} covers the relevant preliminaries, and Section \ref{sec_Main} proves Theorem \ref{thm_Cox_DRP}.  Section \ref{sec_Conclusion} ends with some concluding remarks.

\section{Preliminaries}\label{sec_prelims}

Recall that $S \subseteq [A]^\omega$ is stationary if it meets every closed, unbounded subset of $[A]^\omega$ (in the sense of Jech~\cite{MR1940513}).  By Kueker~\cite{MR0457191} this is equivalent to requiring that for every $f: [A]^{<\omega} \to A$ there is an element of $S$ that is closed under $f$.

In what follows, we will use uppercase letters to denote second order variables/parameters, and lowercase letters to denote first order variables/parameters.  We will also use some standard abbreviations; e.g.\ if our language includes the $\in$ symbol, $v$ is a first order variable, and $Z$ is a second order variable, ``$v = Z$" is short for 
\[
\forall x \ x \in v \ \iff \ Z(x).
\]

Given a structure $\mathfrak{A}=(A,\dots)$ (which we always assume to have a countable signature), the satisfaction relation in Stationary Logic is defined recursively by:
\begin{gather*}
\mathfrak{A} \models \ \text{stat} Z \  \phi(Z,U_1,\dots,U_\ell, p_1,\dots,p_k)  \\
\iff \\
\big\{ Z \in [A]^\omega \ : \ \mathfrak{A} \models \ \phi(Z,U_1,\dots,U_\ell, p_1,\dots,p_k)  \big\} \text{ is stationary in } [A]^\omega.
\end{gather*}

We define a hierarchy of formulas in Stationary Logic that mimics the usual hierarchy in Second Order Logic.  Since 
\[
\text{aa} Z \ \phi(Z,\dots)
\]
roughly translates as
\[
\exists C \ \ C \text{ is club and } \ \forall Z \in C \ \phi(Z,\dots) , 
\]
the \emph{aa} quantifier will correspond to the existential second order quantifier when constructing the hierarchy.  Similarly, since
\[
\text{stat}  Z \ \phi(Z,\dots)
\]
roughly translates as
\[
\forall C \ \ C \text{ is club } \implies \ \exists Z \in C \ \phi(Z,\dots),
\]
the \emph{stat} quantifier will correspond to the universal second-order quantifier.

\begin{definition}
A formula in Stationary Logic without second order quantifiers will be denoted by $\boldsymbol{\Sigma^1_0}$ or $\boldsymbol{\Pi^1_0}$.  For $n > 0$, a formula of the form
\[
\text{stat} Z_1 \ \dots \ \text{stat} Z_k \ \ \phi(Z_1,\dots, Z_k, \dots)
\]
where $\phi$ is $\Sigma^1_{n-1}$ will be called a  $\boldsymbol{\Pi^1_n}$ formula, and a formula of the form
\[
\text{aa} Z_1 \ \dots \ \text{aa} Z_k \ \ \psi(Z_1,\dots, Z_k, \dots)
\]
where $\psi$ is $\Pi^1_{n-1}$ will be called a  $\boldsymbol{\Sigma^1_n}$ formula.
\end{definition}
For example, if $\phi(Z_0,Z_1,v_1,\dots,v_\ell)$ has no \emph{stat} or \emph{aa} quantifers, then
\[
\text{stat} Z_0 \ \ \text{aa} Z_1  \ \ \phi(Z_0,Z_1,v_1,\dots,v_\ell)
\]
is a $\Pi^1_2$ formula.

\begin{definition}\label{def_LST}
We say that the \textbf{LST property holds for Stationary Logic} iff for every structure $\mathfrak{A}=(A,\dots)$\footnote{Recall we always assume countable signature, though for everything discussed in this paper an $\omega_1$-sized signature would still be fine.} there exists a $W \subseteq A$ of size $\le \omega_1$ such that for all formulas $\phi$ in Stationary Logic with no free occurrences of second order variables, and all first order parameters $p_1,\dots,p_k \in W$, 
\[
\mathfrak{A} \models \phi[\vec{p}] \ \text{ if and only if } \mathfrak{A}|W \models \phi[\vec{p}].
\]
We say that the \textbf{$\boldsymbol{\Pi^1_1 \! \! \downarrow}$ LST property holds for Stationary Logic} iff for every structure $\mathfrak{A}=(A,\dots)$ there exists a $W \subseteq A$ of size $\le \omega_1$ such that for all $\Pi^1_1$ formulas $\phi$ in Stationary Logic with no free occurrences of second order variables, and all first order parameters $p_1,\dots,p_k \in W$, 
\[
\textbf{if } \mathfrak{A} \models \phi[\vec{p}] \textbf{, then } \mathfrak{A}|W \models \phi[\vec{p}].
\]
\end{definition}

\begin{remark}
Note that in the definition of the $\Pi^1_1 \!  \! \downarrow$ LST property, we only require that $\Pi^1_1$ formulas reflect \textbf{downward}.  If there is always an $\omega_1$ sized substructure that reflects $\Pi^1_1$ formulas both upward and downward, then the full LST property holds.  This issue is discussed further in Section \ref{sec_Conclusion}.
\end{remark}

We consider variants of the \textbf{Diagonal Reflection Principle} introduced in Cox~\cite{DRP} and \cite{Cox_RP_IS}.  We use the following definition, which by Cox-Fuchs~\cite{Cox_Fuchs_DRP} is equivalent to the definitions from \cite{DRP} and \cite{Cox_RP_IS}:
\begin{definition}\label{def_DRP_internal}
$\text{DRP}_{\text{internal}}$ asserts that for every sufficiently large regular $\theta$, there are stationarily many $W \in \wp_{\omega_2}(H_\theta)$ such that:
\begin{itemize}
 \item $|W|=\omega_1 \subset W$; and
 \item Whenever $A \in W$ is uncountable and $S \in W$ is a stationary subset of $[A]^\omega$, the set $S \cap W \cap [W \cap A]^\omega$ is stationary in $[W \cap A]^\omega$.
\end{itemize}
\end{definition}
The ``internal" part of the definition refers to the fact that we require that $S \cap W \cap [W \cap A]^\omega$ is stationary, not merely that $S \cap [W \cap A]^\omega$ is stationary.  Definition \ref{def_DRP_internal} is simply the diagonal version of an internal variant of WRP introduced in Fuchino-Usuba~\cite{FuchinoUsuba} (see Cox~\cite{Cox_RP_IS} for a discussion).

\section{Proof of Theorem \ref{thm_Cox_DRP}}\label{sec_Main}

We prove a slightly stronger variant of Theorem \ref{thm_Cox_DRP}.  The proof below is strongly influenced by Fuchino et al~\cite{FuchinoEtAl_DRP_LST}.

\begin{theorem}\label{thm_DRP_char}
The following are equivalent:
\begin{enumerate}
 \item\label{item_DRP_internal} $\text{DRP}_{\text{internal}}$.

 \item\label{item_diag_Veryweak} For every structure $\mathfrak{A} = (A,\dots)$, there is a $W \subseteq A$ of size at most $\omega_1$ such that for every finite list $p_1,\dots,p_k \in W \cap A$ and every formula $\phi$ without 2nd order quantifiers,
 \[
\Big( \mathfrak{A} \models \ \text{stat} Z \ \phi[Z,\vec{p}] \Big) \ \implies \ \Big( \mathfrak{A}|W \models \text{stat} Z \ \phi[Z,\vec{p}] \Big).
 \]

 \item\label{item_diag_weak} The $\Pi^1_1 \! \! \! \downarrow$-LST property holds of Stationary Logic (as in Definition \ref{def_LST});

\item\label{item_diag_strong} For every structure $\mathfrak{A}=(A,\dots)$, there is a $W \subseteq A$ of size at most $\omega_1$ such that for every formula $\psi$ in 2nd order prenex form with no free occurrences of second order variables, and every finite list $p_1,\dots, p_k$ of elements of $W$, \textbf{if} 
\[
\mathfrak{A} \models  \psi[\vec{p}]  
\]
then, letting $\hat{\psi}$ be the formula obtained from $\psi$ by changing all \emph{aa} quantifiers to \emph{stat} quantifiers, 
\[
\mathfrak{A}|W \ \models \ \hat{\psi}[\vec{p}].
\]

\end{enumerate}
\end{theorem}

Before proving the theorem, we remark that in parts \ref{item_diag_Veryweak}, \ref{item_diag_weak}, and \ref{item_diag_strong} of Theorem \ref{thm_DRP_char}, we only mentioned first order parameters from $W \cap A$.  If the structure $\mathfrak{A}$ is sufficiently rich then it often makes sense to also speak of second-order parameters that are elements of $W$.  But in general (e.g.\ when $\mathfrak{A}$ is a group) it is more natural to only speak of first order parameters from $W \cap A$. 


\begin{proof}
(of Theorem \ref{thm_DRP_char}):  \eqref{item_diag_strong} trivially implies \eqref{item_diag_weak}, since if $\psi$ is represented as a prenex $\Pi^1_1$ formula, then $\hat{\psi} = \psi$ (because there are no \emph{aa} quantifiers in the original formula at all).  Similarly, \eqref{item_diag_weak} trivially implies \eqref{item_diag_Veryweak} because if $\phi$ has no second order quantifiers, 
\[
\text{stat} Z \ \phi
\]
 is obviously a $\Pi^1_1$ formula.

To see that \eqref{item_diag_Veryweak} implies \eqref{item_DRP_internal}, assume \eqref{item_diag_Veryweak} and suppose $\theta$ is a regular cardinal $\ge \omega_2$.  We need to find a $W \prec (H_\theta,\in)$ such that $|W|=\omega_1 \subset W$ and for every $s \in W$ that is a stationary collection of countable sets,
\[
s\cap W \cap \left[W \cap \bigcup s \right]^\omega \text{ is stationary.}
\]
Consider $\mathfrak{A} = (H_\theta,\in)$.  Let $W \subset H_\theta$ be as in the statement of \eqref{item_diag_Veryweak}.  Fix any $s \in W$ that is a stationary collection of countable sets.  Then
\[
\mathfrak{A} \models \ \text{stat} Z \ \exists p \ p= Z \cap \bigcup s \text{ and } p \in s
\]
and hence, since $s \in W$ and the only second order quantifier in the (prenex) formula above is a \emph{stat} quantifier, 
\[
\mathfrak{A}|W \  \models \ \text{stat} Z \ \exists p \ p= Z  \cap \bigcup s \text{ and } p \in s. 
\]
Unravelling the definition of the satisfaction relation, this means that 
\[
\big\{  Z \in [W]^\omega \ : \   Z \cap \bigcup s   \in W \cap s  \big\} \text{ is stationary in } [W]^\omega
\]
and it follows that $W \cap s \cap \big[ W \cap \bigcup s \big]^\omega$ is stationary in $\big[ W \cap \bigcup s \big]^\omega$.

To see that $\omega_1 \subset W$, it suffices  to show that $W \cap \omega_1$ is uncountable (since by first-order elementarity of $W$ in $(H_\theta,\in)$, $W \cap \omega_1$ is transitive).  Now
\[
\mathfrak{A} \models \text{stat} Z \ \ \exists p \ \exists \alpha \ \ \big( p = Z \cap \omega_1, \ \alpha <  \omega_1, \text{ and } \alpha \text{ is an upper bound of } p\big),
\]
so by assumption on $W$, this statement is also satisfied by $\mathfrak{A}|W$ (note that the parameter $\omega_1$ is an element of $W$ because $\omega_1$ is first-order definable in $\mathfrak{A}$ and $W$ is at least first-order elementary in $\mathfrak{A}$).  If $W \cap \omega_1$ were countable, say $W \cap \omega_1 = \delta < \omega_1$, it would follow that for stationarily many $Z \in W \cap [W]^\omega$, there is an $\alpha < W \cap \omega_1 = \delta$ such that $\alpha$ is an upper bound of $Z \cap \delta$.  This would be a contradiction, since due to the countability of $\delta$, the set of $Z \in [W]^\omega$ such that $\delta \subseteq Z$ is a club.

Finally, to prove that \eqref{item_DRP_internal} implies \eqref{item_diag_strong}:  fix a structure $\mathfrak{A}=(A,\dots)$ and let $\theta$ be a sufficiently large regular cardinal with $\mathfrak{A} \in H_\theta$.  By \eqref{item_DRP_internal} there is a $W \prec (H_\theta,\in,\mathfrak{A})$ witnessing $\text{DRP}_{\text{internal}}$.   We prove by induction on complexity of formulas $\psi$ in 2nd order prenex form that if $p_1,\dots, p_k \in W \cap A$ and
\[
\mathfrak{A} \models \psi[\vec{p}]
\] 
then, letting $\hat{\psi}$ be the result of replacing all \emph{aa} quantifiers with \emph{stat} quantifiers,
\[
\mathfrak{A}|(W \cap A) \models \hat{\psi}[\vec{p}].
\]
We actually need to inductively prove a slightly stronger statement:  namely, that whenever $\psi$ is a 2nd order prenex formula, $p_1,\dots,p_k \in W \cap A$, and $Z_1,\dots,Z_\ell \in W \cap [A]^\omega$,
\begin{equation}
\mathfrak{A} \models \ \psi[\vec{Z},\vec{p}] \ \implies \ \ \mathfrak{A}|(W \cap A) \models \hat{\psi}[\vec{Z},\vec{p}]. 
\end{equation}
So suppose 
\begin{equation}\label{eq_A_sat}
\mathfrak{A} \models \ QZ \ \phi[Z,U_1,\dots,U_k, p_1,\dots,p_\ell]
\end{equation}
where $Q$ is either the \emph{aa} or \emph{stat} quantifier, $U_1,\dots,U_k$ are each elements of $W \cap [A]^\omega$, $p_1,\dots,p_\ell \in W \cap A$, and the inductive hypothesis holds of the formula $\phi$.

Now regardless of whether $Q$ is the \emph{aa} or \emph{stat} quantifier, 
\[
\widehat{QZ \ \phi} \ \equiv \ \text{stat} Z \ \hat{\phi}.
\]
and by \eqref{eq_A_sat}  (since the \emph{aa} quantifier is stronger than the \emph{stat} quantifier) 
\[
\mathfrak{A} \models\ \text{stat}Z \ \phi[Z,U_1,\dots,U_k,p_1,\dots,p_\ell].
\]
Hence, by the definition of the stationary logic satisfaction relation,
\[
s:= \Big\{ Z \in [A]^\omega \ :  \   \mathfrak{A} \models \phi[Z,\vec{U},\vec{p} \ ]  \Big\} \text{ is stationary in } [A]^\omega.
\]
Note that since $\vec{U}$, $\vec{p}$, $\phi$, and $\mathfrak{A}$ are elements of $W$, it follows that $s \in W$.  Since $W$ is internally diagonally reflecting,
\[
s \cap W \cap [W \cap A]^\omega \text{ is stationary in } [W \cap A]^\omega.
\]

Consider for the moment an arbitrary $Z \in s \cap W \cap [W \cap A]^\omega$.  Then
\[
\mathfrak{A} \models \phi[Z,\vec{U},\vec{p}]
\]
and it follows by the induction hypothesis (and that $Z$, $\vec{U}$, and $\vec{p}$ are each elements of $W$) that:
\[
\mathfrak{A}|(W \cap A) \models \hat{\phi}[Z,\vec{U},\vec{p}].
\]

Hence, we have shown that
\[
s \cap W \cap [W \cap A]^\omega \subseteq \big\{ Z \in [W \cap A]^\omega \ : \ \mathfrak{A}|(W \cap A) \models \hat{\phi}[Z,\vec{U},\vec{p}]  \big\}.
\]
Since the set on the left side is stationary, the set on the right side is too.  So by the definition of the satisfaction relation,
\[
\mathfrak{A}|(W \cap A) \models \ \text{stat} Z \ \hat{\phi}[Z,\vec{U},\vec{p} \ ].
\]
This completes the proof of the \eqref{item_DRP_internal} $\implies$ \eqref{item_diag_strong} direction.
 
\end{proof}

\section{Concluding remarks}\label{sec_Conclusion}

We remark that it is straightforward to show, in ZFC alone, that:
\begin{lemma}\label{lem_ZFC_Sigma11}
For every structure $\mathfrak{A}=(A,\dots)$ there exists a $W \subseteq A$ of size at most $\omega_1$ such that 
\[
\mathfrak{A}|W \ \prec^{\Sigma^1_1}_{\downarrow} \ \mathfrak{A}
\]
(i.e.\ such that $\Sigma^1_1$ formulas satisfied by $\mathfrak{A}$ are also satisfied by $\mathfrak{A}|W$).

In fact, if $\theta$ is a regular cardinal such that $\mathfrak{A} \in H_\theta$, and 
\[
W \prec_{\text{1st order}} (H_\theta,\in,\mathfrak{A})
\]
 is such that $|W|=\omega_1$ and
\begin{equation}\label{eq_internallyClub}
W \cap [W \cap A]^\omega \text{ contains a club in } [W \cap A]^\omega
\end{equation}
(this always holds for stationarily many $W$, e.g.\ for those $W$ that are internally approachable), then 
\[
\mathfrak{A}|(W \cap A) \ \prec^{\Sigma^1_1}_{\downarrow} \ \mathfrak{A}.
\]
\end{lemma}
We briefly sketch the proof of the lemma; more details, and other related results, can be found in Cox~\cite{RP_EklofEtAl}.  One proves by induction on complexity of formulas, making use of \eqref{eq_internallyClub}, that if $\phi$ is $\Sigma^1_1$, $p_1,\dots,p_k \in W \cap A$, and $Z_1,\dots, Z_\ell \in W \cap [A]^\omega$, then
\[
\textbf{if } \mathfrak{A} \ \models \phi[\vec{Z},\vec{p}] \ \textbf{, then } \mathfrak{A}|(W \cap A) \models \phi[\vec{Z},\vec{p}].
\]
This was basically part of the proof from Fuchino et al~\cite{FuchinoEtAl_DRP_LST} that $\text{DRP}_{\text{IC}}$ implied the LST for Stationary Logic.  See \cite{RP_EklofEtAl} for some other related ZFC theorems.

So by Lemma \ref{lem_ZFC_Sigma11} one can always get an $\omega_1$ sized substructure that reflects all $\Sigma^1_1$ statements downward.  And if $\text{DRP}_{\text{internal}}$ holds, one can \emph{also} get an $\omega_1$ sized substructure that reflects all $\Pi^1_1$ statements downward.  But it is consistent that both of these are true, yet no \emph{single} $\omega_1$-sized substructure downward reflects all $\Pi^1_1$ \emph{and} all $\Sigma^1_1$ statements.  In particular, in any model where $\text{DRP}_{\text{internal}}$ holds and $\text{DRP}_{\text{IC}}$ fails, Theorem \ref{thm_Cox_DRP} tells us that there is a structure such that no $\omega_1$-sized substructure reflects all $\Pi^1_1$ and all $\Sigma^1_1$ statements (though there are structures that reflect one or the other).  

Another way to view this phenomenon, in terms of DRP-like principles, is that $\text{DRP}_{\text{internal}}$ yields stationarily many $W \in \wp_{\omega_2}(H_\theta)$ such that the transitive collapse $H_W$ of $W$ is ``correct about stationary sets"; i.e.\ whenever $s \in H_W$ and $H_W \models$ ``$s$ is a stationary set of countable sets", then $V$ believes this too. However, if $W$ is not internally club, it is possible (by \cite{Cox_RP_IS}) that $H_W$ is correct about stationary sets, but is \emph{not} correct about clubs; i.e.\ there can be a $c \in H_W$ such that $H_W \models$ ``$c$ is a club of countable sets", but $V$ does not believe this.  If, on the other hand, $W$ witnesses $\text{DRP}_{\text{IC}}$, then $H_W $ is correct about \emph{both} stationarity \emph{and} clubness.

\end{document}